\documentclass{birkjour}

\usepackage[all]{xy}

\usepackage{graphicx}
\usepackage{amsthm}
\usepackage{amstext}
\usepackage{amsmath,amscd,amsfonts}
\usepackage{amssymb}
\usepackage{mathrsfs}
\usepackage{mathtools}
\usepackage[all]{xy}
\usepackage{color}
\usepackage[table]{xcolor}
\usepackage{comment}
\usepackage{enumerate}
\newcommand{\ran}{\operatorname{ran}}
\newcommand{\alg}{\operatorname{alg}}
\newcommand{\co}{\operatorname{co}}

\newcommand{\A}{\mbox{${\mathcal A}$}}

\newcommand{\E}{\mbox{${\mathcal E}$}}
\newcommand{\F}{\mbox{${\mathcal F}$}}

\newcommand{\I}{\mbox{${\mathcal I}$}}
\newcommand{\J}{\mbox{${\mathcal J}$}}

\newcommand{\Ll}{\mbox{${\mathcal L}$}}

\newcommand{\N}{\mbox{${\mathcal N}$}}
\renewcommand{\O}{\mbox{${\mathcal O}$}}
\renewcommand{\P}{\mbox{${\mathcal P}$}}
\newcommand{\Q}{\mbox{${\mathcal Q}$}}

\renewcommand{\S}{\mbox{${\mathcal S}$}}
\newcommand{\T}{\mbox{${\mathcal T}$}}
\newcommand{\U}{\mbox{${\mathcal U}$}}

\newcommand{\intern}[1]{\langle #1\rangle}
\newcommand{\sspan}{\operatorname{span}}

\newcommand{\naturais}{\mathbb {N}}

\newtheorem{theorem}{Theorem}[section]
\newtheorem{lemma}[theorem]{Lemma}
\newtheorem{corollary}[theorem]{Corollary}
\newtheorem{proposition}[theorem]{Proposition}
\theoremstyle{definition}

\theoremstyle{remark}
\newtheorem{remark}[theorem]{Remark}
\newtheorem{example}{Example}

\numberwithin{equation}{section}

\def\sqr#1#2{{\,\vcenter{\vbox{\hrule height.#2pt\hbox{\vrule width.#2pt
height#1pt \kern#1pt\vrule width.#2pt}\hrule height.#2pt}}\,}}

\begin{document}

\title[On a class of left ideals of nest algebras]{On a class of left ideals of nest algebras}

\author[P. Costa]{Pedro Costa}
\address{Mathematical Institute, University of Oxford\\
Andrew Wiles Building, Radcliffe Observatory Quarter, Woodstock Road, Oxford OX2 6GG, United Kingdom}
\email{Pedro.BezerraRoqueDaCosta@maths.ox.ac.uk}

\author[M. Ferreira]{Martim Ferreira}
\address{
Instituto Superior T\'ecnico \\
Universidade de Lisboa\\
Av. Rovisco Pais 1\\
1049-001 Lisbon, Portugal}
\email{martim.ferreira@tecnico.ulisboa.pt}

\author[L. Oliveira]{Lina Oliveira}
\address{Centro de An\'alise Matem\'atica, Geometria e Sistemas Din\^amicos\\
Departamento de Matem\'atica\\
Instituto Superior T\'ecnico \\
Universidade de Lisboa\\
Av. Rovisco Pais 1\\
1049-001 Lisbon, Portugal}
\email{lina.oliveira@tecnico.ulisboa.pt}

\begin{abstract} 
We introduce a class of left ideals (and subalgebras) of nest algebras determined by totally ordered families of partial isometries on a   
complex Hilbert space $H$. 
 Let $\E$ be a  family of partial isometries that is totally ordered in the Halmos--McLaughlin ordering, and let 
$\A_{\E}$ be the subset of operators  in $B(H)$ which, for all $E\in \E$, map the initial space of $E$ to the final space of $E$. 
We show that  $\A_{\E}$   is a  subalgebra of $B(H)$ if and only if $\A_{\E}$ is  a  left ideal of a  certain  nest algebra, and if so,    
 $\E$ consists  of  power partial isometries, except possibly for its supremum $\vee \E$, in which case  the range $\ran(\vee\E)$ is $H$. 
It is also shown that any  left ideal $\A_{\E}$  is decomposable  and  that the subset of finite rank operators in  its closed unit ball  is strongly dense in the ball. Necessary and sufficient conditions to solve  $Tx=y$ and $T^*x=y$ in $\A_{\E}$ are given.
\end{abstract}

\thanks{The first author was supported by a Calouste Gulbenkian Foundation grant under the \emph{Novos Talentos em Matem\'atica} programme and ``laCaixa'' LCF/BQ/EU24/12060096. The second author was  funded by FCT/Portugal through project UIDB/04459/2020 with DOI identifier 10-54499/UIDP/04459/2020. The third author was partially supported  by FCT/Portugal through project UIDB/04459/2020 with DOI identifier 10-54499/UIDP/04459/2020.}



\maketitle

\section{Introduction}\label{s_intro}
Ideals and, more generally, bimodules of nest algebras have attracted a considerable attention over the years, with their characterisation  being achieved for the associative, Jordan, and Lie structures (cf. \cite{Cap, DO, EP, HMS, LiLi}). Until quite recently, the emphasis of the research seems to have been given essentially to  two-sided ideals and  (bi)modules, being   their one-sided versions  not that much present in the literature. However, in a very recent  paper on the ideal structure of  nest algebras, J.L. Orr  established  the form of  left ideals of a nest algebra $\alg \N$ associated with a nest $\N$ of projections  on a separable complex Hilbert space $H$, by identifying them with the so-called constructible left ideals. That is, a left ideal  $\Ll\subseteq \alg \N$ is said to be \emph{constructible} if there exists a net $(N_\alpha, x_\alpha)_{\alpha\in A}$, indexed by some set $A$, where 
$N_\alpha\in \N$ and $x_\alpha\in H$ are such that 
$$
\Ll=\{T\in \alg \N\colon \lim_{\alpha\in A} \|(I-N_\alpha)Tx_\alpha\|=0\}
$$
(see \cite{Orr}).
 
 The present work furthers the investigation into the structure of left ideals and, at the same time, introduces   a new class of operator algebras said $2$-nest algebras. Without assuming the separability of the Hilbert space, we provide a way of constructing a subclass  of  left ideals (and algebras), notwithstanding  the approach  being different from that of \cite{Orr}: here the left ideals are characterised by means of operators leaving ``invariant'' or ``preserving'' a totally ordered family of partial isometries in a sense to be made clear  (see \eqref{e3} below). These   left ideals of a nest algebra are classified in a way reminiscent of ideals and bimodules of nest algebras (see Section 2), allowing for a simple and intuitive way of obtaining concrete examples, which we will do in the sequel. The key strategy to obtain  this description of left ideals is moving from projections to partial isometries, as will be seen.

 Let $H$ be a complex Hilbert space and let $B(H)$ be the set of bounded linear operators on $H$. Recall that an operator $V\in B(H)$ is said to be a \emph{partial isometry} if $V$ is an isometry on the orthogonal complement $\ker(V)^{\perp}$ of its kernel $\ker(V)$. That is to say that, for all $x\in \ker(V)^\perp$, 
 $\lVert Vx \rVert = \lVert x \rVert.$

The set   $\U(H)$ of partial isometries on the Hilbert space $H$ is endowed with a partial ordering. This partial order was  introduced by P.R. Halmos and J.E. McLaughlin in \cite{HM}  and is defined as follows.
 Given partial isometries $E,F\in \U(H)$, $E$ \emph{is less than or equal to} $F$, denoted  $E\leq F$,  if, for all $x$ in the initial space of $E$, we have  $F(x)=E(x)$. It is clear from this definition   that the initial (respectively, final) space of $E$ is included in the initial (respectively, final) space of $F$,  both inclusions being strict whenever $E\neq F$.
This partial ordering extends that of  self-adjoint projections on $H$ but, unlike   projections,  partial isometries do not form a lattice. In fact, a pair of partial isometries might fail to have an upper bound.  

The poset of partial isometries has been shown to have far-reaching properties and importance in the structure of  $B(H)$. For example, L. Moln\'ar  proved that bijective maps on $\U(H)$ (with $\dim H\geq 3$) which preserve order in both directions can be extended to $B(H)$ as linear or conjugate-linear maps, under some other not too strict  conditions (see  \cite{Molnar}). 

More recently, a  motivation for considering bijections preserving order in the context of Physics was given by Y. Friedman and A.M. Peralta. In particular, the authors proved  that a bijection preserving the ordering  of partial isometries in both directions (with some other mild constraints) could be extended to $B(H)$ as a real-linear Jordan triple isomorphism. In fact, the setting in which the results were obtained is of a  more general nature than  $B(H)$ (see  \cite{FP}).

The set $\U(H)$ of partial isometries can  be endowed with other partial orders (see, for example,  \cite{Ross}), however, here, we restrict ourselves to the Halmos--McLaughlin ordering. 

This work is divided into four sections that we now describe.

 Notwithstanding our main focus be on complete  totally ordered families of partial isometries,  in Section 1 we begin by considering general posets of partial isometries. A totally ordered family $\E$ of partial isometries  in $\U(H)$ is said to be \emph{complete} if it contains $0$ and arbitrary infima and suprema.
 
 A bounded linear operator $T\in B(H)$ is said to \emph{preserve} an isometry $E\in \U(H) $ if  $TE^*E(H)\subseteq EE^*(H)$, that is, if $T$ maps the initial space of $E$ to the final space of $E$.  If $T$ preserves all the partial isometries of $\E$, then we say that $T$ \emph{preserves} $\E$.
 
 In Section 2, Theorem \ref{t_algebra}, we show that the   space $\A_{\E}$ consisting of all operators in $B(H)$ preserving a complete totally ordered family $\E$ of partial isometries is an algebra if and only if it is a left ideal of a  nest algebra determined by $\E$. In fact, under these circumstances, $\E$ itself must consist  of power partial isometries, except possibly for its supremum $\vee \E$, in which case the range $\ran(\vee\E)$ is $H$. Recall that power partial isometries, that is, partial isometries $V$ whose positive powers $V^n$ are also  partial isometries, were completely characterised by P.R. Halmos and L.J. Wallen in \cite{HW}: 
 
 \begin{theorem}[\cite{HW}, Theorem]\label{HaWa}
Every power partial isometry is a direct sum whose summands are unitary operators, pure isometries, pure co-isometries, and truncated shifts. The direct sum representation can be expressed in such a way that each type of summand (and, in particular, each index for a truncated shift) occurs at most once; the representation so expressed is unique.
\end{theorem} 
 
 A subalgebra $\A_{\E}$ of $B(H)$ associated with  some totally ordered family $\E$   of  partial isometries
is said to be  a \emph{$2$-nest algebra}. These weakly closed algebras, although having always an ``{\sl enveloping}'' nest algebra, do properly generalise nest algebras in the infinite dimensional setting (see Sections 2 and 4).  Nest algebras and   $2$-nest algebras  do  have an  intrinsically different structure (see Example \ref{exa2}) but, 
 for  finite dimensional Hilbert spaces, nest algebras and $2$-nest algebras coincide (see Remark \ref{rem3}). Nest algebras and $2$-nest algebras differ exactly when the latter are non-unital (see Corollary \ref{c_identity}).
 
 The set of finite rank operators is investigated in Section 3. It is proved that every rank $n$ operator in a $2$-nest algebra $\A_{\E}$ is the sum of $n$ rank one operators in $\A_{\E}$, in other words, $\A_{\E}$ is decomposable. Furthermore, the subset of finite rank operators in the closed unit ball of $\A_{\E}$ is strongly dense in the ball. 
 
 Section 4 is devoted to analysing the solvability of  $Tx=y$ and $T^*x=y$ in a $2$-nest algebra, for which we provide necessary and sufficient conditions.
  
  We give several examples illustrating the theory throughout the manuscript. On  occasion, when referring to  others' results, we opted for their transcription, to ease the readability of the present work.
  
  In the remainder of this section, we establish the notation and a few facts concerning the Halmos--McLaughlin ordering of partial isometries. 
  
In what follows, $H$ is always a finite or infinite dimensional complex Hilbert space, $B(H)$ is the set of bounded linear operators on $H$, and $\U(H)$ is the poset of partial isometries on $H$. Given a \emph{nest} in $B(H)$, that is, a totally ordered family $\N$ of (self-adjoint) projections  which contains $0$ and the identity $I$ on $H$, the corresponding \emph{nest algebra} $\alg\N$ consists of the operators in $B(H)$ leaving invariant the ranges of all projections in $\N$ (see \cite{Ringrose}). For a comprehensive account of the theory of nest algebras, the reader is referred to \cite{Davidson}.

The next lemma is a well-known characterisation of partial isometries (see \cite{Halmos2}, for example), of which we make a note here for later reference.

\begin{lemma}\label{l_basic}
Let $V \in \mathcal{B(H)}$. The following assertions  are equivalent.
\begin{enumerate}[(i)]
    \item $V$ is a partial isometry.
    \item $V^*$ is a partial isometry.
    \item $VV^*$ is a projection.
    \item $V^*V$ is a projection.
    \item $VV^*V = V$.
    \item $V^*VV^* = V^*$.
\end{enumerate}
\end{lemma}
It immediately follows  from (v) and (vi) of this lemma that $V^*V, VV^*$ are self-adjoint (equivalently, orthogonal)  projections
on the (closed) subspaces $\ker(V)^{\perp}$ and $\ran(V)$, respectively. 
The subspaces $\ker(V)^{\perp}$   and $\ran(V)$ are called, respectively,  the \emph{initial space} and the  \emph{final space} of $V$.    The set of partial isometries  on $H$ will be denoted by  $\U(H)$. 

For a partial isometry $V$, denote by $P_V$ the \emph{initial projection} $V^*V$ onto the initial space of  $V$, and denote by $Q_V$  the  \emph{final projection} $VV^*$ onto the final space of  $V$.

\begin{proposition}\label{p_bounds} Let $\F$ be  a subset  of partial isometries in $\U(H)$. The following hold.
\begin{enumerate}[(i)]
\item $\F$ has an infimum $\wedge \F$ in $\U(H)$.
\item If $\F$ is bounded above, then $\F$ has a supremum $\vee \F$ in $\U(H)$.
\item If $\F$ is totally ordered then $\F$ has an infimum and a supremum in $\U(H)$.
\end{enumerate}
\end{proposition}

\begin{proof} (i) 
Let $X$ be the closed subspace defined by
    $$X = \Bigl\{  x\in\bigcap\limits_{\substack{E\in \F}} \ker(E)^{\perp} \colon E(x)= F(x), \, \mbox{for all}\,  E, F \in \F \Bigr\}.$$
Clearly, $X$ is the largest closed subspace of $H$ on which all the partial isometries in $\F$ agree. In consequence, the bounded linear operator $T$ defined, for all $x\in H$, by 
    $$
    Tx=\begin{cases}
			Ex, & x \in X\\
            0, & x \in X^{\perp}
		 \end{cases},
    $$
    where $E\in \F$ is a (any) fixed partial isometry is a partial isometry, and $T$ is the largest partial isometry below every element of $\F$. Hence $T=\wedge \F$.
    
    (ii)    Let $U$ be an upper bound for $\F$. Since $U$ lies above each partial isometry $E$  in $\F$, it follows that the set $X$ defined by  
    $$X = \overline{\sspan} \Bigl\{\bigcup\limits_{\substack{E\in \F}} \ker(E)^{\perp} \Bigr\}$$
     is contained in $\ker(U)^{\perp}$. 
      It is clear that the partial isometry $V$ defined by 
     $$
    Vx=\begin{cases}
			Ux, & x \in X\\
            0, & x \in X^{\perp}
		 \end{cases}
    $$
is the supremum of $\F$.
 
 (iii) By (i) and (ii), it suffices to show that $\F$ has an upper bound.
 Notice that, since $\F$ is totally ordered, any two partial isometries in $\F$ agree on the initial space of the smallest. 
 
 Let $X=\overline{ \bigcup\limits_{\substack{E\in \F}} \ker(E)^{\perp} }$, and define the operator $T$ on $H=X\oplus X^\perp$ as follows. If $x\in X^\perp$, then set $T(x)=0$; if $x\in X$, then 
 $ T(x)=\lim E_n(x_n)$, where $(x_n)$ is a (any) sequence in $X$ converging to $x$ in norm, and the sequence $(E_n)$ is such that, for all $n\in \naturais$, the element $x_n$ lies in the initial space of $E_n$.  
 By continuity, $T$ is an isometry whose  initial space is $X$. Moreover, by construction, $T$ agrees with any $E\in \F$ on the initial space of the latter. Hence $T$ is an upper bound (in fact, the supremum) of $\F$, concluding the proof.
 \end{proof}

A family $\mathcal{F}$  
 of partial isometries is said to be a \emph{complete family } if $0\in \mathcal{F}$ and it contains all arbitrary infima and suprema, that is, 
given any subset   of $\mathcal{F}$, its supremum and infimum exist and lie in  $\mathcal{F}$.

\begin{lemma}\label{l_co}
Let $\F$ be a totally ordered family of partial isometries on a complex Hilbert space $H$. Then,  there exists a smallest totally ordered complete family of partial isometries containing  $\F$. Moreover, this family of partial isometries consists of $\{0\}$ and all the  infima and suprema of the subsets of $\F$.\end{lemma}

\begin{proof}
Let $\S$ be the class  of all totally ordered families of partial isometries containing $\mathcal{F}$ and $\{0\}$. The class $\S$ can be partially ordered by inclusion, and every chain $\{ \mathcal{F}_i \}_{i\in I}$ of $\S$ is bounded by $\bigcup\limits_{\substack{i\in I}} \mathcal{F}_i \in \S$. Therefore, by Zorn's lemma, there exists a maximal element of $\S$, i.e., $\mathcal{F}_1$ such that $\mathcal{F}\subseteq \mathcal{F}_1$ and $\mathcal{F}_1$ is maximal.

Furthermore, $\mathcal{F}_1$ must be complete. Otherwise, we could add  a missing infimum or supremum of a subset of $\F_1$, ending up with another  totally ordered family containing $\F$ and strictly containing $\F_1$, which would contradict the maximality of $\F_1$. 

Since the intersection of totally ordered complete families is also a totally ordered complete family, it follows that there exists a smallest totally ordered complete family containing $\mathcal{F}$. 

The remaining assertion is an immediate consequence of the definition of complete family and of the reasoning above.
\end{proof}

Given a totally ordered family $\F$ of partial isometries, the \emph{complete cover} $\co(\mathcal{F})$ of $\F$ is defined to be the smallest complete family of partial isometries containing $\F$.

\begin{remark}\label{r_3}
Notice that, when $\F$ is totally ordered, so are the corresponding families $\{P_E\colon E\in \F\}$ of initial projections and  $\{Q_E\colon E\in \F\}$ final projections. 
Consequently, the sets $\N_P=\{P_E\colon E\in \F\}\cup\{0, I\}$ and $\N_Q=\{Q_E\colon E\in \F\}\cup \{0, I\}$ are nests of projections in $B(H)$. 
In what follows, $\alg \N_P$ and $\alg \N_Q$ will denote the corresponding nest algebras.
\end{remark}

\section{Operators preserving partial isometries}\label{s_oppi}
A bounded linear operator $T\in B(H)$ is said to \emph{preserve} a partial isometry $E\in \U(H) $ if  $TE^*E(H)\subseteq EE^*(H)$, that is, if $T$ maps the initial space of $E$ to the final space of $E$. Equivalently,  
 $$TE^*E=EE^*TE^*E \iff(I-Q_E)TP_E=0.
 $$
 If $T$ preserves all the partial isometries of a 
given  
 family $\F$, then we say that $T$ \emph{preserves} $\F$. 
 
 Let  $\A_{\F}$ be the subset of $B(H)$ consisting of the operators  that preserve $\F$:
\begin{equation}\label{e3}
    \mathcal{A}_{\F} = \{ T \in \mathcal{B(H)} :  (I-Q_E)TP_E = 0, \forall E\in \mathcal{F} \}.
\end{equation}

Notice that if $\F$ is totally ordered, then  $\F$ itself is contained in ${\A}_{\F}$. In fact, given  a fixed $E\in \F$, we have, for all $F\in \F$, 

$$
(I-Q_F)EP_F=(I-Q_F)F=0, \quad \mbox {if}\; F\leq E;
$$
and
$$
(I-Q_F)EP_F=(I-Q_F)E=(I-Q_F)Q_EE=0, \quad \mbox {if}\; E\leq F,
$$
since in this case $I-Q_F\leq I-Q_E$. Hence, $E$ lies in $\A_{\F}$.

\begin{lemma}\label{l_cover}
Let $\mathcal{F}$ be a totally ordered family of partial isometries on a complex Hilbert space $H$. Then $\mathcal{A}_{\mathcal{F}} = \mathcal{A}_{\co(\mathcal{F})}$.
\end{lemma} 
\begin{proof}
Since  $\mathcal{F}\subseteq \co(\mathcal{F})$,  we have $\mathcal{A}_{\co(\mathcal{F})} \subseteq \mathcal{A}_{\mathcal{F}}$. 
As to the reverse inclusion, notice that, by Lemma \ref{l_co}, it suffices to show that $T$ preserves the infimum and supremum of arbitrary subsets of $\F$. This is immediate, by the continuity of $T$.
\end{proof}

 The set $\A_{\F}$ is a complex vector subspace of $B(H)$ but not necessarily a subalgebra, as the next simple example shows.
\begin{example}\label{ex1}
Let $\mathcal{F} =  \{V\}$ where $V\colon \mathbb{C}^{2} \to \mathbb{C}^{2}$ is the partial isometry   such that  $(z,w)\mapsto (0,z)$. Let  $T\colon \mathbb{C}^{2} \to \mathbb{C}^{2}$ 
be the operator defined, for all  $(z,w)\in \mathbb{C}^{2}$, by $T(z,w)=(w,z)$. 
The operator $T$ lies in $\mathcal{A}$.  However, $TV \notin \mathcal{A}$, since $TV(1,0) = T(0,1) = (1,0)$  does not lie in the final space of $V$, yielding that $\mathcal{A}$ is not an algebra.
\end{example}

 We have nevertheless necessary and sufficient conditions on $\F$ for $\A_{\F}$ to be an algebra.
\begin{theorem}\label{t_algebra}
Let $H$ be a complex Hilbert space, let $\mathcal{F}$ be a totally ordered family of partial isometries in $B(H)$, and let $\A_{\F}$ be as in \eqref{e3}. The following are equivalent.
\begin{enumerate}[(i)]
    \item $\mathcal{A}_{\F}$ is an algebra.
    
 \item $\mathcal{A}_{\F}$ is a left ideal of the nest algebra $\alg \N_Q$.
    \item For all $E\in \mathcal{F}$, $\ran(E) \subseteq \ker(E)^{\perp}$ or  $\ran(E) = H$.
 \end{enumerate}
\end{theorem}

It is easily seen that $\A_{\F}$ is a weakly closed subspace, that  is, $\A_{\F}$ is closed in the weak operator topology. Notice also that, if there exists $E\in \F$ for which $\ran(E) = H$, then $E=\vee \F$.

Before giving the proof of the theorem, we need  the following  lemma.

\begin{lemma}\label{l_xy}Let $\F$ be a  
  totally ordered family of partial isometries on a complex Hilbert space $H$. 
Let $E\in\mathcal{F}$ be a partial isometry that simultaneously satisfies $\ran(E)\neq H$ 
 and $\ran(E)\not\subseteq \ker(E)^{\perp}$. 
 Then, there exist  vectors $x,y \in \mathcal{H}$ such that $x \notin \ran(E)^{\perp}$, $y\notin \ran(E)$, and, for all $V\in \mathcal{F}$, $x\in\ker(V)$ or $y\in \ran(V)$.
\end{lemma}

\begin{proof} Since $\ran(E)\nsubseteq \ker(E)^{\perp}$, there exists  $z\in \ran(E)$ such that $z = x + z_{\perp}$, with $z_{\perp}\in\ker(E)^{\perp}$, $x\in\ker(E)$,  and $x\neq 0$. 
Observe that this implies that $x\notin \ran(E)^{\perp}$, since $\intern{z,x}=\|x\|^2\neq 0$. 

Consider the (totally ordered) subset  $\S_x = \{V\in\mathcal{F} \colon x\in \ker(V)\}$ of $\F$. Observe that $E$ lies in $\S_x$. We analyse separately the two possible situations:  $\S_x=\F$, and  $\S_x\subsetneq \F$.

If $\S_x=\F$, it suffices to take any non-zero $y\in \ran(E)^{\perp}$, to conclude the proof.

Suppose now that  $\S_x\subsetneq \F$, and  
let $W$ be the supremum of $\S_x$. 
 If $E<W$, then any $y\in \ran(W)\ominus \ran(E)$ satisfies the conditions. 
 The only other possibility is $E=W$. In these circumstances, $x\in \ker(V)$ if and only if $V\leq E$.

Case  a) There exists $V\in \F$, with $E<V$,  such that $\ran(E) \nsubseteq \ker(V)^{\perp}$. 
 
If this is the case, then there exists  $z'\in\ran(E)$ such that $z' = x' + z'_{\perp}$, with $x'\in\ker(V)$, $z'_{\perp}\in\ker(V)^{\perp}$, and $x'\neq 0$. 
  Observe that $x' \notin \ran(E)^{\perp}$, since $\intern{z', x'} = \| x' \|^2 \neq 0$. 
 
 As above, define $S_{x'}$. Since $V \in S_{x'}$, either $S_{x'} = \F$, and it suffices to take any non-zero $y \in \ran(E)^{\perp}$, or $\vee S_{x'} \geq V > E$, and it suffices to choose any $y \in \ran( \vee S_{x'} ) \ominus \ran(E)$.

Case  b) For all $V\in\mathcal{F}$,  with $E<V$, we have $\ran(E)\subseteq\ker(V)^{\perp}$. Then,   
$\ran(E)\subseteq \bigcap\limits_{\substack{E<V}}\ker(V)^{\perp}$. Since 
$\ran(E)\nsubseteq \ker(E)^{\perp}$, there exists a non-zero closed subspace $X\subseteq \bigcap\limits_{\substack{E<V}}\ker(V)^{\perp}$, such that 
$$
\bigcap\limits_{\substack{E<V}}\ker(V)^{\perp} = \ker(E)^{\perp} \oplus X.$$

Observe that each $V > E$  is an isometry (thus injective) on $\bigcap\limits_{\substack{E<V}}\ker(V)^{\perp}$ and maps $\ker(E)^{\perp}$ onto $\ran (E)$. It follows that, it suffices to choose a non-zero vector $w\in X$ and set $y = Vw$, for any $V>E$, since all $V$ agree on $w$.
\end{proof}

Next, we  prove  Theorem \ref{t_algebra}.
\begin{proof} 
$(i) \Longrightarrow(ii)$: Suppose that $T$ is an operator in $\A_{\F}$, and let $E$ be a partial isometry in $\F$. We want to show that, for all $y\in H$, 
$$Q_ETQ_Ey=TQ_Ey.$$

Given $y\in H$, there exists $x\in H$ such that $Q_Ey=Ex$. Then, Lemma \ref{l_basic}~(v)
$$
Q_ETQ_Ey = Q_ETEx = Q_ETEP_Ex = TEP_Ex 
    = TEx.
$$
Hence, 
$$
Q_ETQ_Ey = TEx=TQ_Ey.
$$
The space  $\mathcal{A}_{\mathcal{F}}$ is a left ideal of $ \alg \N_Q$, since, given $E\in \E$, for $T\in\mathcal{A}_{\mathcal{E}},S\in\alg \N_Q$, we have  $ST(\ker E^{\perp}) \subseteq S(\ran E) \subseteq \ran E$.

$(ii) \Longrightarrow(iii)$:
Suppose that there exists a partial isometry $E\in\mathcal{F}$ that simultaneously satisfies $\ran(E)\neq H$ and $\ran(E)\not\subseteq \ker(E)^{\perp}$.  By Lemma \ref{l_xy}, there exist non-zero vectors $x,y \in \mathcal{H}$ such that $x \notin \ran(E)^{\perp}$, $y\notin \ran(E)$ and, for all $V\in \mathcal{F}$, $x\in\ker(V)$ or $y\in \ran(V)$.
 
 The rank one operator $T = x\otimes y$ lies in $\A_{\F}$.   
 In fact,  for any $V\in\mathcal{F}$, if $x\in\ker(V)$, $T(\ker(V)^{\perp})=0$, otherwise $y\in\ran(V)$ and $T(\mathcal{H})\subseteq \ran(V)$. Either way, $T(\ker(V)^{\perp})\subseteq \ran(V)$.
However, since $x \notin \ran(E)^{\perp}$ and $y\notin \ran(E)$, it follows that $T(\ran(E)) = \sspan\{ y\}\nsubseteq \ran(E)$. Hence,  $T\notin \alg \N_Q$.

$(iii) \Longrightarrow(i)$:
Assume that, for all $E\in \mathcal{F}$, we have $\ran(E) \subseteq \ker(E)^{\perp}$ or $\ran(E) = H$. 
 We show next that, for  $A,B \in \mathcal{A}_{\F}$, 
 we have that $AB$ lies also in $\A_{\F}$. Notice that, for $E\in \F$ with $\ran(E) = H$, we have  immediately that, $ AB( \ker(E)^{\perp})\subseteq  \ran(E)$. Consequently, it only remains to consider the case when  $E\in\F$ is such that $\ran(E) \subseteq \ker(E)^{\perp}$. 
 
 Let $x\in \ker(E)^{\perp}$. Then, $Bx\in \ran(E)$ from which follows that  $Bx \in \ker(E)^{\perp}$. Hence, $ ABx = A(Bx)\in \ran(E)$, which concludes the proof that $\A_{\F}$ is closed under the product on $B(H)$.
\end{proof}

Recall that a partial isometry $V$ in $B(H)$ is said to be a \emph{power partial isometry} if, for all $n\in \naturais$, its powers $V^n$ are also partial isometries (see \cite{HW}). 

An immediate consequence of Theorem \ref{t_algebra} is the following corollary.

\begin{corollary}\label{c_algebra}
Let $H$ be a complex Hilbert space, let $\mathcal{F}$ be a complete totally ordered family of partial isometries in $B(H)$, and let $\A_{\F}$ be as in \eqref{e3}. If $\mathcal{A}_{\F}$ is an algebra, then   (i) the family $\F$ consists entirely of power partial isometries, or (ii) all partial isometries in $\F$ are power partial isometries except for $\vee\F$ and $\ran(\vee\F)=H$.
 \end{corollary}

\begin{remark}\label{rem3}
If $H$ is  finite dimensional, then for any partial isometry $E\in B(H)$, we have, by the rank-nullity theorem,  that $\dim \ker(E)^\perp=\dim \ran(E)$. Hence, if $\ran(E)\subseteq \ker(E)^\perp$, then these two spaces must coincide. If $\ran E=H$, then $E$ is surjective, hence bijective, and again $\ran(E)\subseteq \ker(E)^\perp$.  As a consequence of Theorem \ref{t_algebra}, we have that, given a  totally ordered family $\F$ of partial isometries on a finite dimensional Hilbert space $H$, the operator space $\A_{\F}$ is  an algebra if and only if it is  a nest algebra. 

However, this is not the case, if $H$ is infinite dimensional. For example, let $H=\ell_2(\naturais)$ be the space of square summable complex sequences, and let $E=S_f^2P$ be the partial isometry where $S_f$ is the forward shift operator and $P$ is the projection on $\{(x_n)\in \ell_2(\naturais)\colon x_{2n+1}=0, \; \mbox{for all}\; n\in \naturais\}$. If  $\F=\{E\}$, then, by Theorem \ref{t_algebra}, the space $\A_{\F}$ is an algebra but it is not a nest algebra, since it does not contain the identity operator. 
\end{remark}

\begin{remark} \label{remextra}
Corollary \ref{c_algebra} cannot be improved inasmuch as a family $\F$ of power partial isometries does not guarantee that $\A_{\F}$ is an algebra. For example, 
 consider the complete totally ordered family $\F$ of partial isometries $0 \leq ... \leq V_n \leq ... \leq V_1$ in  $H=\ell_2(\mathbb{Z})$, where  
$$ V_1 (..., x_{-4}, x_{-3}, x_{-2}) = (..., x_{-5}, x_{-4}, x_{-3}),$$
    $$ V_1(x_{-1}, x_0, x_1) = (x_0, x_1, x_{-1}), $$
    $$ V_1(x_2, x_3, x_4, ...) = (0, x_2, x_3,...), $$  
  and, for all $j \neq 1$, define $V_j(e_i) = 0$, whenever $2 \leq i \leq j$, and $V_j(e_i) = V_1(e_i)$, otherwise. 
  It follows from Theorem \ref{HaWa} that $V_1$ is a power partial isometry, as are the remaining operators in $\F$.
     However, the space $\A_{\F}$ is not an algebra, as it does not satisfy  Theorem  \ref{t_algebra}~(iii).

\end{remark}

 Since, by Theorem \ref{t_algebra}, the space  $\A_{\F}$ is an algebra exactly when it is a left ideal of some nest algebra, one has  in fact that $\A_{\F}$ is a nest algebra if and only if it contains the identity operator.

\begin{corollary}\label{c_identity}
Let $H$ be a complex Hilbert space, let $\mathcal{F}$ be a complete totally ordered family of partial isometries in $B(H)$, and let $\A_{\F}$ be as in \eqref{e3}. Then, $\A_{\F}$ is a nest algebra if and only if it contains the identity  $I$ on $H$.
\end{corollary}

In what follows, we will be concerned with the set $\A_{\F}$ associated with a family $\F$ of partial isometries as in \eqref{e3}. Lemma \ref{l_cover} allows for considering only complete families $\F$, which we shall do henceforth. The families of partial isometries that are complete and totally ordered will be denoted by $\E$.

An algebra $\A_{\E}$ determined by a (complete) totally ordered family  $\E$ of partial isometries as in \eqref{e3} will be called a \emph{$2$-nest algebra}.

\begin{example}\label{exa8} Let $H=L^2([0,1])$, and consider the totally ordered family $\E=\{\Phi_t\colon t \in [0,1]\}$ of partial isometries,  where for all $t \in [0,1]$, the partial isometry $\Phi_t$ is defined by

    \[
        \Phi_t(f) = 
            \begin{cases}
                0 & \text{ if} \int_0^t \lvert f \rvert^2 = 0, \\
                \phi(f) & \text{ otherwise},
            \end{cases}
    \]

    with 
     \[
        \phi(f(x)) = 
            \begin{cases}
            2\phi(f(2x)) & \text{ if\, } x\in [0,\tfrac12],\\
                0 &  \text{ otherwise}.
            \end{cases}
    \]
 It is clear that  $\Phi_0 = 0$.

    We have, for all $t\in ]0,1]$, that
$$
P_{\Phi_t} (H)= \ker{\Phi_t}^\perp = \{ f \in L^2([0,1]): f(x) = 0 \text{ a.e. in } [t,1] \} =: N_t 
$$
contains  
$$
 Q_{\Phi_t}(H) = \ran{\Phi_t} = N_{t/2} 
 $$
Hence, by Theorem \ref{t_algebra}, the space $\A_{\E}$ is a $2$-nest algebra.
\end{example}
\begin{example}
 Consider the family of partial isometries $0 \leq ... \leq V_n \leq ... \leq V_1$ in $\ell^2(\mathbb{Z})$ defined  as follows:
$$ V_1 (..., x_{-4}, x_{-3}, x_{-2}) = (..., x_{-4}, x_{-3}, x_{-2}, 0); $$
    $$ V_1(x_{-1}, x_0, x_1) = (x_0, x_1, x_{-1}); $$
    $$ V_1(x_2, x_3, x_4, ...) = (0, x_2, x_3,...) .$$ 
 Now, for $j \neq 1$,  define   $V_j(e_i) = 0$, for all  $2 \leq i \leq j$ and $V_j(e_i) = V_1(e_i)$, otherwise.  Notice that, by Theorem \ref{HaWa}, the operator $V_1$ is a  power partial isometry and, hence, so are the remaining operators.
    Applying  Theorem  \ref{t_algebra}, we have now that $\A_{\E}$ is a $2$-nest algebra.
    \end{example}

\section{Finite rank operators}\label{s_finiterank}
As is well-known, the finite rank operators in a nest algebra are decomposable and are dense in the strong operator topology. We will show next that the same is true for 2-nest algebras. The proofs of Lemma \ref{ef} and Theorem \ref{t_decomp} are inspired by those existing for nest algebras (see \cite{Erdos}).

Firstly, we fix some notation. Given a partial isometry $E$ in a totally ordered complete family of partial isometries $\mathcal{E}$, define 
$$
E_- = \vee\{V\in \E\colon  V< E\}, \qquad E_+= \wedge\{V\in\mathcal{E}\colon  E< V\}.
$$

\begin{lemma}\label{ef}
Let $H$ be a complex Hilbert space, let $\mathcal{E}$ be a totally ordered complete family of partial isometries in $B(H)$,  let $U= \vee\mathcal{E}$, and let $\mathcal{A}_{\mathcal{E}}$ be the set of operators in \eqref{e3}.
Then, given   
vectors $e,f\in {H}$,
 the 
  operator $e\otimes f$ lies in $\mathcal{A}_{\mathcal{E}}$ if and only if 
$e\in\ker(U)$, or there exists an $E\in\mathcal{E}$ such that $f\in\ran(E)$ and $e\in \ker(E_- )$.
\end{lemma}
\begin{proof} The lemma holds trivially when $e=0$ or $f=0$. Hence, suppose that $e,f\neq 0$.

Let $e\otimes f$ lie in $\A_{\E}$, and assume that $e\notin\ker(U)$. Then, since $e\otimes f\in\mathcal{A}_{\E}$, we have that $f\in\ran(U)$. Notice that, if it were the case that $e\in\ker(U)$, then $e\otimes f$ would coincide with zero on the initial space of each partial isometry in $\E$.

Let  $E = \wedge\{V\in\mathcal{F}\colon  f\in\ran(V)\}$. Notice that $E\neq 0$, since $f\in \ran(E)$.

Let $V\in\mathcal{F}$ be such that $V<E$. Since $f\notin\ran(V)$ and $e\otimes f (\ker(V)^{\perp})\subseteq \ran(V)$, it follows that,  for all $x\in\ker(V)^{\perp}$, we must have $\langle x,e \rangle = 0$. In other words, $e\in\ker(V)$, for any $V\in \E$ lying below $E$. 
Hence,   for each $V\in\mathcal{F}$ such that $V<E$, we have that $e$ is orthogonal to $\ker(V)^{\perp}$ and, thus, also to $\ker(E_-)^{\perp}$.

The converse is immediate.
\end{proof}

We make here a note of a lemma needed in the proof of Theorem \ref{t_decomp} below.

\begin{lemma}\cite[Lemma 1]{Erdos}\label{4:3}
Let $\N$ be  a  
 nest in a complex  Hilbert space ${H}$, and let $x$ be a fixed vector in ${H}$. Then, the mapping $E\mapsto Ex$ from $\N$ to ${H}$ is continuous with respect to the order topology of $\N$.
\end{lemma}

We are now able to prove the decomposability of the finite rank operators in a 2-nest algebra.

\begin{theorem}\label{t_decomp} Let $H$ be a complex Hilbert space,  let $\mathcal{E}$ be a totally ordered complete family of partial isometries in $B(H)$, and let 
 $\mathcal{A}_{\mathcal{E}}$ be a $2$-nest algebra. If  $R$ is  an operator of rank $n\in\naturais$ in $\mathcal{A}_{\mathcal{E}}$, then $R$ can be written as the sum of $n$ rank one operators in  $\mathcal{A}_{\mathcal{E}}$.
\end{theorem}
\begin{proof}
Let  
$R = \sum_{i=1}^nx_i\otimes y_i$
 be a rank $n$ operator in $\mathcal{A}_{\mathcal{E}}$, where $\{x_i \colon 1\leq i\leq n\}$ and $\{y_i \colon 1\leq i\leq n\}$ are linearly independent sets. 

Given a partial isometry $E \in \mathcal{E}$, we have,  for all $x\in {H}$,
\begin{equation}\label{e6}
0=(I-EE^*)RE^*Ex=
\sum_{i=1}^n\langle x, E^*Ex_i \rangle(I-EE^*)y_i.
\end{equation}

Let $\P$ and $\Q$ be the sets of initial and final spaces of the partial isometries in $\mathcal{E}$, respectively. Notice that $\P$ and $\Q$ are totally ordered and complete  lattices, that are homeomorphic to $\mathcal{E}$ in the order topology.

Let $\mathcal{E}_1$ be the subset of $\mathcal{E}$ consisting of the partial isometries $E$ such that $E^*Ex_i=0$, for $1\leq i\leq n$. Let $\mathcal{E}_2$ be the subset of $\E$ consisting of the partial isometries $E$ for which 
the set  $\{(I-EE^*)y_i \colon 1\leq i\leq n\}$ is linearly dependent. Since $\{y_i\colon 1\leq i\leq n \}$ is a linearly independent set, we have that $\mathcal{E}_2\neq \mathcal{E}$.  

It is clear from \eqref{e6} that $\mathcal{E}_1 \cup \mathcal{E}_2 = \mathcal{E}$. 
 Moreover, by  Lemmas \ref{4:3}, continuity, and the fact that $\P$ and $\mathcal{E}$ are homeomorphic,  we have that $\mathcal{E}_1$ is closed in the order topology. 
 
  Recall that a finite set  of vectors is linearly dependent if and only if its
Grammian determinant  is zero. Hence, by   Lemmas \ref{4:3}, continuity, and the fact that $\Q$ and $\mathcal{E}$ are homeomorphic,  it follows that $\mathcal{E}_2$ is closed.

Suppose firstly that  $\mathcal{E}_2 = \emptyset$. Hence, $\mathcal{E}_1 = \mathcal{E}$, and $\{x_i\colon 1\leq i\leq n\}\subseteq \ker( \vee_{V\in\mathcal{E}} V )$ . It follows from Lemma \ref{ef} that, for all $1\leq i\leq n$, the rank one operators $x_i\otimes y_i$ lie in $ \mathcal{A}_{\mathcal{E}}$, as desired.

Assume now that $\mathcal{E}_2\neq \emptyset$, and let $U = \wedge\{E \colon E\in\mathcal{E}_2\}$. Since $\mathcal{E}_2$ is closed, $U\in\mathcal{E}_2$ and $\{(I-UU^*)y_i\colon  1\leq i\leq n\}$ is a linearly dependent set of vectors. Suppose, without loss of generality, that
$$
(I-UU^*)y_1 = \sum\limits_{i=2}^n\alpha_i(I-UU^*)y_i,
$$
where the $\alpha_i$'s are complex numbers. Using this, we have

\begin{align*}
    R&= x_1 \otimes \bigl[ UU^*y_1 + (I-UU^*)y_1 \bigr]+\sum\limits_{i=2}^n x_i\otimes y_i\\
    &= x_1 \otimes z + \sum\limits_{i=2}^N (x_i + \overline{\alpha_i}x_1)\otimes y_i,
\end{align*}
where
$$z = UU^*y_1 - \sum\limits_{i=2}^n\alpha_iUU^*y_i.$$

Since $\mathcal{E} = \mathcal{E}_1 \cup \mathcal{E}_2$ and $\mathcal{E}_1$ is closed, it is easy to see from the definition of $U$ that $U_- = \vee\{ V\in\mathcal{E}\colon V<U\}$ lies in $\mathcal{E}_1$. Consequently $(U_-)^*U_-x_1 = 0$, that is, $x_1\in\ker(U_-)$. 

On the other hand, 
 since $z\in\ran(U)$, it follows from Lemma \ref{ef} that $x_1\otimes z \in \mathcal{A}_{\mathcal{E}}$ . Hence,  $R$ can be written as a sum of a rank one operator in $\mathcal{A}_{\mathcal{E}}$ and a rank $n-1$ operator in $\mathcal{A}_{\mathcal{E}}$. Now, the  proof can be  completed using a simple induction argument.
\end{proof}

\begin{theorem}
Let $\mathcal{A}_{\E}$ be a $2$-nest algebra associated with a totally ordered complete family $\E$ of partial isometries,  let $\mathcal{R}_1$ be the set of finite rank operators in the  closed unit ball of $\A_{\E}$. Then, the SOT-closure of $\mathcal{R}_1$  coincides with the closed unit ball of $\mathcal{A}_{\E}$.
\end{theorem}
\begin{proof}
Let $I$ denote the identity operator on $H$. Since $\alg \N_Q$ is a nest algebra and $I\in\alg N_Q$, by \cite[Theorem 3]{Erdos}, there exists a net $(R_j)$ of finite rank operators in the unit ball of $\alg \N_Q$ which converges  to $I$ in the strong operator topology.  Hence, given any $T$ in the unit ball of $\mathcal{A}_{\mathcal{E}}$, by Theorem \ref{t_algebra}~(ii), the net $(R_jT)$ consists of finite rank operators in the unit ball of $\mathcal{A}_{\mathcal{E}}$ and converges  to $T$ in the strong operator topology, as required.
\end{proof}

\section{The equations  $Tx=y$ and $T^*x=y$}
In this section, we address the possibility of solving the equations $Tx=y$  and $T^*x=y$ in a $2$-nest subalgebra of $B(H)$, for some fixed $x,y\in H$. We will  extend the results and techniques developed in \cite{Lance} to obtain necessary and sufficient conditions to solve said equations  in a $2$-nest algebra. 
However, before doing that, we will analyse next an example which motivates  what follows and, at the same time, gives a glimpse of how $2$-nest algebras and nest algebras may  differ. 

\begin{example}\label{exa2} Let $0 \leq ... \leq E_n \leq ... \leq E_1$ be the   totally ordered complete family $\E$ of partial isometries in $\ell_2(\naturais)$ defined, for all $j\in \naturais,$ by
$$
E_j(e_i) = 
            \begin{cases}
                e_{2i} & \text{if } i>j \\
                0 & \text{if } i\leq j
            \end{cases},
$$
where $(e_i)_{i\in \naturais}$ is the standard Hilbert basis of $\ell_2(\naturais)$.
By Theorem \ref{t_algebra}, the operator space $\A_{\E}$ is a $2$-nest algebra, since, for all $j\in \naturais$, 
\begin{equation*}
        \ran(Q_{E_j}) = \ran({E_j}) = \overline{\sspan\{ e_{2i}\colon i > j\}}\subseteq \overline{\sspan\{ e_i\colon i > j\}}= \ker(E_j)^\perp .
    \end{equation*}
    
    It is easily seen that the  set $\I=\{T\in \A_{\E} \colon \intern{Te_1, e_1}=0\}$  is a weakly closed ideal of the $2$-nest algebra $\A_{\E}$.   
    
    However, one cannot ``mimic'' the characterisation of weakly closed ideals  of nest algebras, as $\I$ cannot be defined by means of an  order homomorphism  on $\E$ along the lines of \cite{EP}. 
    In fact, if one considers an order homomorphism $E\mapsto \tilde{E}$ on $\E$ and the space 
    $\mathcal{J}=\{T\in \A_{\E}\colon (I - Q_{\tilde{E}}) T P_E = 0\}$, it is immediately seen that, given $T\in\J$, it does not follow necessarily that $\intern{Te_1,e_1}=0$.     For  details on ideals and, more generally, bimodules of nest algebras and their connection with order homomorphisms, see \cite{EP}. 
\end{example}

A constructive description   of the order homomorphism on the nest associated with a given weakly closed bimodule was given  in \cite{Cap} in the Banach space setting. It was of utmost relevance in this construction how the rank one operators $e\otimes f$ in the bimodule behaved under left (and right) multiplication by an operator $T$ in the nest algebra. That is, if $T$ lies in a nest algebra $\A$ and $e\otimes f$ were to be in a $\A$-bimodule, then 
$e\otimes Tf$ and $T^*e\otimes f$ must also lie in the bimodule (see \cite{Cap}, Lemma 4 and Theorem 1). 
 More generally, what can we say about  the  equations $Tx=y$ and $T^*x=y$? When are they solvable? 
The remainder of this section is concerned with these questions.

The proof of the next lemma can be found in \cite{Lance}. The space of the $n\times n$ upper triangular complex matrices is denoted by $\T_n$.

\begin{lemma} \label{l_eqs1}
Let $x = \left( x_1, x_2, \dots, x_n \right)$ and $y = \left( y_1, y_2, \dots, y_n \right)$ be vectors in $\mathbb{C}^n$, and suppose that there exists $K \geq 0$ such that, for all integers   $1 \leq r \leq n$, 
       $ \sum_{i = r}^n \lvert y_i \rvert^2 \leq K^2 \sum_{i = r}^n \lvert x_i \rvert^2.$ 
    Then, there exists a matrix $A$ in $\mathcal{T}_n$ such that $\lVert A \rVert \leq K$ and $Ax = y$. Moreover, 
    \begin{itemize}
	\item [(i)] if  $y_p = 0$, for some  $1 \leq p \leq n$, then the  matrix $A = (a_{ij})\in \T_n$ can be chosen such that $a_{pj} = 0$, for all $1 \leq j \leq n$;
    
        \item [(ii)] if $x_q = 0$, for some $1 \leq q \leq n$, then $A = (a_{ij})\in \T_n$ can be chosen such that $a_{iq} = 0$, for $1 \leq i \leq n$.
    \end{itemize}

\end{lemma}

\begin{lemma}\label{l_operator} Let $\E$ be a complete totally ordered family of partial isometries on a complex Hilbert space $H$.
    Let $n$ be a  positive integer, and let $\E_n = \{E_i: i\in \naturais, \, 1 \leq i \leq n \}$ be a  complete finite subset of $\mathcal{E}$, and let the mapping $i\mapsto E_i$ be injective and order preserving. Let $x \in P_{E_n}, y \in Q_{E_n}$, and let $\mathcal{A}_{\E_n}$ be the operator space associated with $\E_n$. Suppose that there exists $K \geq 0$ such that, for all $E \in {\E_n}$, 
    \begin{equation}\label{10}
        \lVert (I - Q_E) y \rVert \leq K \lVert (I - P_E) x \rVert.
    \end{equation}
    Then, there exists $T_n \in \mathcal{A}_{\E_n}$ such that $T_n x = y$ and $\lVert T_n \rVert \leq K$.
    \end{lemma} 
   
\begin{proof}
    For $1 \leq i \leq n$, define
$$x_i = P_{E_i} (I - P_{E_{i-1}}) x = P_{E_i}x - P_{E_{i-1}}x$$
$$y_i = Q_{E_i} (I - Q_{E_{i-1}}) y = Q_{E_i}y - Q_{E_{i-1}}y$$

and $H_x = \text{span } \{ x_i, 1 \leq i \leq n \}$.
It is easily seen that  $\{x_i: 1 \leq i \leq n\}$ is an orthogonal subset of $H$ as is $\{y_i: 1 \leq i \leq n\}$.

It is also the case that,   for all $1 \leq i,j \leq n$ with $i\leq j$, the operator $x_j \otimes y_i \in \mathcal{A}_{\E_n}$. Indeed, given $1 \leq k \leq n$, we have
\begin{equation*}
        (I - Q_{E_k}) (x_j \otimes y_i) P_{E_k} =P_{E_k} x_j \otimes (I - Q_{E_k}) y_i ,
\end{equation*}
from which follows that, if $k \geq i$,
$$
(I - Q_{E_k}) y_i =
             (I - Q_{E_k}) Q_{E_i} (I - Q_{E_{i-1}}) y = 
             (Q_{E_i} - Q_{E_i})  (I - Q_{E_{i-1}}) y = 0 
$$
and, if $k < i$, 
$$
  P_{E_k} x_j = 
             P_{E_k} P_{E_j} (I - P_{E_{j-1}}) x = 
             P_{E_k} (I - P_{E_{j-1}}) x = 
             (P_{E_k} - P_{E_{k}}) x = 0.
$$
Hence, for all $1 \leq i,j \leq n$ with $i\leq j$, the operator $x_j \otimes y_i \in \mathcal{A}_{\E_n}$.

    Let $A = (a_{ij}) \in \mathcal{T}_n$ be an upper triangular matrix, and let $T_A$ be the operator on $H$  defined by 
    $$
    T_A = \sum_{1 \leq i \leq j \leq n} a_{ij} (x_j \otimes y_i).
    $$
    Since $T_A$ is a linear combination of operators  in $\mathcal{A}_{\E_n}$, it follows that 
   $T_A$ lies in  $\mathcal{A}_{\E_n}$. Moreover,   since $T_A$ coincides with zero on $H_x^\perp$, we have  $\lVert T_A \rVert = \sup_{\substack{z \in H_x \\ \lVert z \rVert = 1}} \lVert T_A z \rVert$.
Setting $
        z = \sum_{i = 1}^n \gamma_i x_i$, we have 
 
\begin{equation*}
        \begin{split}
            \lVert T_A \rVert^2 & = \sup_{\substack{z \in H_x \\ \lVert z \rVert = 1}} \lVert \sum_{1 \leq i \leq j \leq n} a_{ij} \langle z, x_j \rangle y_i \rVert^2  
             = \sup_{\substack{z \in H_x \\ \lVert z \rVert = 1}} \lVert \sum_{1 \leq i \leq j \leq n} a_{ij} \gamma_j \lVert x_j \rVert^2 y_i \rVert^2  \\
              &=\sup_{\substack{z \in H_x \\ \lVert z \rVert = 1}}   \sum_{i,j=1}^n \bigl| a_{ij} \gamma_j \lVert x_j \rVert^2 \bigr|^2 \lVert y_i \rVert^2   
          = \sup_{\substack{z \in H_x \\ \lVert z \rVert = 1}}   \sum_{i,j=1}^n \bigl| \left( a_{ij} \lVert x_j \rVert \lVert y_i \rVert \right) \gamma_j \lVert x_j \rVert \bigr|^2.
           \end{split}
    \end{equation*}
 Consequently,         
       $$
        \lVert T_A \rVert^2  \leq 
             \sup_{\substack{w \in \mathbb{C}^n \\ \lVert w \rVert = 1}}   \sum_{i,j=1}^n \lvert\left( a_{ij} \lVert x_j \rVert \lVert y_i \rVert \right) w_j \rvert^2.
       $$

Let   $B = (b_{ij})$ be the $n\times n$ upper triangular matrix defined by $b_{ij} = a_{ij} \lVert x_j \rVert \lVert y_i \rVert$. The previous calculations show that $\lVert T_A \rVert \leq \lVert B \rVert$.

We show next that $\left( \lVert x_1 \rVert, \lVert x_2 \rVert, \dots, \lVert x_n \rVert \right)$ and $\left( \lVert y_1 \rVert, \lVert y_2 \rVert, \dots, \lVert y_n \rVert \right)$ satisfy the conditions of Lemma \ref{l_eqs1}.
 Let $r$ such that $1 \leq r \leq n$. Then,
$$
        \sum_{i = r}^n \lVert y_i \rVert^2  = \lVert \sum_{i=r}^n y_i \rVert^2  
             = \lVert \sum_{i=r}^n Q_{E_i} (I - Q_{E_{i-1}}) y \rVert^2
             = \lVert (Q_{E_n} - Q_{E_{r-1}}) y \rVert^2.
$$
Similarly, we have that $\sum_{i = r}^n \lVert x_i \rVert^2 = \lVert (P_{E_n} - P_{E_{r-1}}) x \rVert^2$.

 Notice that $y \in Q_{E_n}$ and $x \in P_{E_n}$ and, consequently, 
 $$\sum_{i = r}^n \lVert y_i \rVert^2 = \lVert (I - Q_{E_{r-1}}) y \rVert^2, \qquad \sum_{i = r}^n \lVert x_i \rVert^2 = \lVert (I - P_{E_{r-1}}) x \rVert^2.$$ 
    Hence, it follows from \eqref{10} that 
    $$\sum_{i = r}^n \lVert y_i \rVert^2 \leq K^2 \sum_{i = r}^n \lVert x_i \rVert^2.
    $$
Hence, by Lemma \ref{l_eqs1},   there exists $B \in \mathcal{T}_n$ such that, for all $i=1,\dots, n$,   we have  \begin{equation}\label{11}
\sum_{j=i}^n b_{ij} \lVert x_j \rVert = \lVert y_i \rVert,
\end{equation}
with $\lVert B \rVert \leq K$. 
Observe that Lemma \ref{l_eqs1} also allows for assuming that $b_{ij} = 0$, whenever $x_j=0$ or $y_i=0$.

Bearing the above considerations in mind, we may define an $n\times n$ upper triangular matrix $A = (a_{ij})$ by
\begin{equation*}
        a_{ij} = \begin{cases}
                    \frac{b_{ij}}{\lVert x_j \rVert \lVert y_i \rVert}, & \text{when $b_{ij} \neq 0$} \\
                    0, & \text{when $b_{ij} = 0$}
                \end{cases}
    \end{equation*}
  Hence,   the finite rank operator $T_A$ lies in $\mathcal{A}_{\E_n}$ and $\lVert T_A \rVert \leq \lVert B \rVert \leq \lVert K \rVert$.
    
    Equality \eqref{11} yields, for all $i=1,\dots, n$, 
    $$\sum_{j=i}^n a_{ij} \lVert x_j \rVert^2 \lVert y_i \rVert = \lVert y_i \rVert,
    $$
    from which follows that
     $\sum_{j=i}^n a_{ij} \lVert x_j \rVert^2 = 1$, whenever $\lVert y_i \rVert \neq 0$. Consequently, 
     $$T_A x = \sum_{1 \leq i \leq j \leq n} a_{ij} \langle x, x_j \rangle y_i = \sum_{1 \leq i \leq j \leq n} a_{ij} \lVert x_j \rVert^2 y_i = \sum_{i=1}^n y_i = y.
     $$
The proof is concluded by setting $T_n=T_A$.
\end{proof}
\begin{theorem} \label{t_tx=y}
Let $\E$ be a complete totally ordered family of  partial isometries on a complex Hilbert space $H$, and let $\A_{\E}$ be the operator space corresponding to $\E$.
 Let $x, y \in H$ be such that there exist  $E,F \in \mathcal{E}$ with $ x \in P_E, y \in Q_{F}$. 
 
 The following assertions are equivalent.
    
    \begin{itemize}
    
        \item [(i)] There exists an operator $T \in \mathcal{A}_\mathcal{E}$ such that $Tx = y$.
        
        \item [(ii)] \begin{equation}\label{21}K = \sup_{E \in \mathcal{E}} \frac{\lVert (I - Q_E)y \rVert}{\lVert (I - P_E)x \rVert} < \infty.
        \end{equation}
    \end{itemize}

    Under these conditions, $T$ can be chosen with $\lVert T \rVert = K$. Furthermore, if $S \in \mathcal{A}_\mathcal{E}$ is such that $Sx = y$, then $\lVert S \rVert \geq K$.
\end{theorem} 
    
  It is understood here that the constant $K$ in \eqref{21} is finite, if $(I - Q_E)y=0=(I - P_E)x$, and infinite, if $(I - Q_E)y\neq 0$, $(I - P_E)x=0$.
  
\begin{proof} 

(i) $\Rightarrow$ (ii) Suppose that there exists $S \in \mathcal{A}_\mathcal{E}$ for which $Sx=y$. Then,
for all $E\in \E$,
    \begin{equation}\label{12}
            \lVert (I - Q_E)y \rVert  = \lVert (I-Q_E)Sx \rVert 
           = \lVert (I-Q_E)S (I-P_E) x \rVert
             \leq \lVert S \rVert \lVert (I-P_E) x \rVert,
    \end{equation}
from which follows that $K< \infty.$

(ii) $\Rightarrow$ (i)
 Let  $\Phi$ be the class of all finite complete subfamilies of $\mathcal{E}$ such that, for all $\phi \in \Phi$ we have that  $y \in Q_{\vee\phi}, x \in P_{\vee\phi}$, where $\vee\phi$ denotes the supremum of the finite subfamily $\phi$. 
 By Lemma \ref{l_operator}, for each $\phi \in \Phi$,  there exists  an operator $T_\phi$ such that $T_\phi x = y$,  with $\lVert T \rVert \leq K$. Since the closed ball of radius $K$ and centre $0$ is  WOT-compact, it follows that the net $(T_\phi)$ has an accumulation point $T$ in the weak operator topology, in this ball. Hence, for all $z\in H$, 
$$
\intern{ y,z}=\intern{T_\phi x,z} \rightarrow \intern{Tx,z},
$$
where $\|T\|\leq K$. Hence, $T(x)=y$.

Fix $\phi_0\in \Phi$. Observe that, for $\phi_0, \phi\in \Phi$ with $\phi_0\subseteq \phi$, we have that $\mathcal{A}_\phi \subseteq \mathcal{A}_{\phi_0}$. Since $\mathcal{A}_{\phi_0}$ is weakly closed, it follows  that  $T_\phi \in \mathcal{A}_{\phi_0}$. 

Hence, since $\phi_0 \in \Phi$ is arbitrary, one has that $T \in \cap_{\phi \in \Phi} \mathcal{A}_\phi = \mathcal{A}$, as required.

We have shown that (i)$\iff$(ii). By \eqref{12}, we know that $\|T\|\geq K$, hence $\|T\|=K$. It also follows from the  same equality that the remaining assertion holds.
  \end{proof}

  We turn our attention now to the adjoint operator.
  Notice that, if $T\in \A_{\E}$, then 
  $
  EE^*TE^*E=TE^*E, 
  $
   yielding
   $
   E^*ET^*EE^*=E^*ET^*.
   $ 
   Equivalently, 
   \begin{equation}\label{1721}
   (I-P_E)^\perp T^* Q_E^\perp=(P_E^\perp)^\perp T^* Q_E^\perp=0.
   \end{equation}

   Having in mind \eqref{1721}, we obtain Lemma \ref{l_t*} and Theorem \ref{t_t*x=y} below, which are, respectively, the counterparts of Lemma \ref{l_operator} and Theorem \ref{t_tx=y}  for the adjoint operator. 
   
   \begin{lemma}\label{l_t*}
 Let $\E$ be a complete totally ordered family of partial isometries on a complex Hilbert space $H$.
    Let $n$ be a  positive integer, and let $\E_n = \{E_i: i\in \naturais, \, 1 \leq i \leq n \}$ be a  complete finite subset of $\mathcal{E}$, and let the mapping $i\mapsto E_i$ be injective and order preserving. Let $x \in H, y \in P_{E_n}$, and let $\mathcal{A}_{\E_n}$ be the operator space associated with $\E_n$. Suppose that there exists $K \geq 0$ such that, for all $E \in {\E_n}$, 
    \begin{equation}\label{13}
        \lVert P_E y \rVert \leq K \lVert Q_E x \rVert.
    \end{equation}
    Then, there exists $T_n \in \mathcal{A}_{\E_n}$ such that $T^*_n x = y$ and $\lVert T_n \rVert \leq K$.
    \end{lemma}

\begin{proof}
   Let $x\in H, y \in P_{E_n}$ and define, for $1\leq i\leq n$, 

    $$x_i= (Q_{E_{n-i+1}} - Q_{E_{n-i}})x, \qquad 
     y_i = (P_{E_{n-i+1}} - P_{E_{n-i}})y,
    $$
    and 
    $H_x = \text{span } \{ x_i\colon 1 \leq i \leq n \}$. Notice that, for $i\neq j$, we have that $x_i, x_j$ (respectively, $y_i, y_j$) are orthogonal.

It is also the case that, for $i \leq j$, the operator $y_i \otimes x_j$ lies in  $\mathcal{A}_{\E_n}$, since 
for $1\leq k\leq n$, 
$$(I - Q_{E_k}) (y_i \otimes x_j) P_{E_k} = P_{E_k}y_i \otimes (I - Q_{E_k})x_j,$$
$(I - Q_{E_k}) x_j=0$, if $k > n-j$, and $P_{E_k} y_i =0$, if $k \leq n-i$.

    Let $A = (a_{ij}) \in \mathcal{\T}_n$ be an upper triangular matrix, and define $T_A \in \mathcal{A}_{\E_n}$ by
    $$T_A = \sum_{1 \leq i \leq j \leq n} \overline{a_{ij}} y_i \otimes x_j.
    $$
     Hence, $T^*_A = \sum_{1 \leq i \leq j \leq n} a_{ij} x_j \otimes y_i$ and 
     $\lVert T^*_A \rVert = \sup_{z \in H_x} \frac{\lVert T^*_A z \rVert}{\lVert z \rVert}$.
  For $z \in H_x$, consider the decomposition 
  $
        z = \sum_{i = 1}^n \gamma_i x_i. 
   $
    Then,
\begin{equation*}
    \begin{split}
        \lVert T^*_A \rVert^2 & = \sup_{\substack{z \in H_x \\\lVert z \rVert = 1}} \bigl\Vert \sum_{1 \leq i \leq j \leq n} a_{ij} \gamma_j \lVert x_j \rVert^2 y_i \bigr\Vert^2
        = \sup_{\substack{z \in H_x \\ \lVert z \rVert = 1}} \sum_{i=1}^n \bigl\vert  \sum_{j=i}^n a_{ij} \gamma_j \lVert x_j \rVert^2 \bigr\vert^2 \lVert y_i \rVert^2 \\
        & = \sup_{\substack{z \in H_x \\ \lVert z \rVert = 1}} \sum_{i=1}^n \left\lvert  \sum_{j=i}^n \left( a_{ij} \lVert x_j \rVert \lVert y_i \rVert \right) \gamma_j \lVert x_j \rVert \right\rvert^2.
    \end{split}
\end{equation*}          
It follows that
$$
            \lVert T^*_A \rVert^2 
            \leq \sup_{\substack{w \in \mathbb{C}^n \\ \lVert w \rVert = 1}} \sum_{i=1}^n \left\lvert  \sum_{j=i}^n \left( a_{ij} \lVert x_j \rVert \lVert y_i \rVert \right) w_j \right\rvert^2.
$$
      
    Let $B = (b_{ij})$ be the $n\times n$ upper triangular matrix where $b_{ij} = a_{ij} \lVert x_j \rVert \lVert y_i \rVert$. 
We have   $\lVert T^*_A \rVert \leq \lVert B \rVert$.

   We see now that  $\left( \lVert x_1 \rVert, \lVert x_2 \rVert, \dots, \lVert x_n \rVert \right)$ and $\left( \lVert y_1 \rVert, \lVert y_2 \rVert, \dots, \lVert y_n \rVert \right)$ satisfy the conditions of Lemma \ref{l_eqs1}.
 Given $1 \leq r \leq n$, we have that
 $$
            \sum_{i = r}^n \lVert y_i \rVert^2  = \lVert \sum_{i=r}^n y_i \rVert^2 
       = \lVert \sum_{i=r}^n (P_{E_{n-i+1}} - P_{E_{n-i}}) y \rVert^2
             = \lVert (P_{E_{n-r+1}}) y \rVert^2.
$$
Similarly,  $\sum_{i = r}^n \lVert x_i \rVert^2 = \lVert (Q_{E_{n-r+1}}) x \rVert^2$. Hence, by \eqref{13}, we have
$$\sum_{i = r}^n \lVert y_i \rVert^2 \leq K^2 \sum_{i = r}^n \lVert x_i \rVert^2.
$$

 Applying Lemma \ref{l_eqs1}  to $\left( \lVert x_1 \rVert, \lVert x_2 \rVert, \dots, \lVert x_n \rVert \right)$, $\left( \lVert y_1 \rVert, \lVert y_2 \rVert, \dots, \lVert y_n \rVert \right)$, it follows that there exists a matrix $B \in \mathcal{\T}_n$ such that 
 \begin{equation}\label{35}
 \sum_{j=i}^n b_{ij} \lVert x_j \rVert = \lVert y_i \rVert
 \end{equation}
  and $ \lVert B \rVert \leq K$. Moreover, one may also assume that $b_{ij} = 0$, whenever $x_j=0$ or $y_i=0$.

Define a matrix $A = (a_{ij})$ in $\T_n$ by 
 \begin{equation*}
        a_{ij} = \begin{cases}
                    \frac{b_{ij}}{\lVert x_j \rVert \lVert y_i \rVert}, & \text{when $b_{ij} \neq 0$} \\
                    0, & \text{when $b_{ij} = 0$}.
                \end{cases}
    \end{equation*}
Observe that it was  shown above that  $T_A \in \mathcal{A}_{\E}$ and $\lVert T_A \rVert \leq \lVert B \rVert \leq \lVert K \rVert$.

By \eqref{35}, we have $\sum_{j=i}^n a_{ij} \lVert x_j \rVert^2 \lVert y_i \rVert = \lVert y_i \rVert$ yielding  $\sum_{j=i}^n a_{ij} \lVert x_j \rVert^2 = 1$, whenever $\lVert y_i \rVert \neq 0$. Consequently,

$$T^*_A x = \sum_{1 \leq i \leq j \leq n} a_{ij} \langle x, x_j \rangle y_i =\sum_{1 \leq i \leq j \leq n} a_{ij} \lVert x_j \rVert^2 y_i = \sum_{i=1}^n y_i = P_{E_{n}} y = y.
$$
We conclude the proof by setting $T_n=T_A$.
\end{proof}

   \begin{theorem} \label{t_t*x=y}Let $\E$ be a complete totally ordered family of  partial isometries on a complex Hilbert space $H$, and let $\A_{\E}$ be the operator space corresponding to $\E$.
      Let $x, y \in H$ be such that there exist  $E \in \mathcal{E}$ with $y \in P_E$. 
     The following assertions are equivalent.
     
         \begin{itemize}
        \item[(i)] There exists an operator $T \in \mathcal{A}_\mathcal{E}$ such that $T^*x = y$. 
        
       \item [(ii)]  \begin{equation}\label{367}K = \sup_{E \in \mathcal{E}} \frac{\lVert P_E y \rVert}{\lVert Q_E x \rVert} < \infty.
   \end{equation}
        
    \end{itemize}

    Under these conditions, $T$ can be chosen with $\lVert T \rVert = K$. Furthermore, if $S \in \mathcal{A}_\mathcal{E}$ is such that $S^*x = y$, then $\lVert S \rVert \geq K$.
\end{theorem}

  It is understood here that the constant $K$ in \eqref{367} is finite, if $P_Ey=0=Q_Ex$, and infinite, if $(I - Q_E)y\neq 0$, $(I - P_E)x=0$.

      \begin{proof}  (i) $\Rightarrow$ (ii) Suppose that there is $S \in \mathcal{A}_\mathcal{E}$ such that  $S^*x=y$. Then,
\begin{equation}\label{376}
            \lVert P_E y \rVert  = \lVert P_E S^*x \rVert  
             = \lVert P_E S^* Q_E x \rVert  
             \leq \lVert S \rVert \lVert Q_E x \rVert
\end{equation}
and, hence, $K<\infty$.  

(ii) $\Rightarrow$ (i) Consider the class $\Phi$ of all finite complete subfamilies of $\mathcal{E}$ such that, for all $\phi \in \Phi$, we have $y \in P_{\vee \phi}$. By Lemma \ref{l_t*},  for each $\phi \in \Phi$,  there exists $T_\phi$ such that $T^*_\phi x = y$ with $\lVert T \rVert \leq K$. Observing that the closed ball of radius $K$ and centred
at $0$ is WOT-compact,  a reasoning similar to that of the  proof of Theorem \ref{t_tx=y} will yield the existence of $T$ in said ball such that    $T^*x = y$.

 The  remaining assertions follow easily using a reasoning similar to that of the proof of   Theorem \ref{t_tx=y}.
     \end{proof}

  \vspace{3ex}
  
\subsection*{Acknowledgements} The authors wish to thank   Professor Janko Bra\v{c}i\v{c} for valuable discussions during the preparation of the manuscript and Professor Marek Ptak for bringing insightful examples to their attention.

\end{document}